\documentstyle[amssymb]{article}
\let\mathcal\cal
\begin{document}
\title{A general geometric construction for affine surface area}
\author{Elisabeth Werner\thanks{Supported by a grant from the National
Science Foundation.
\newline
\newline
MSC classification:   52
}}

\date{ }
\maketitle

\begin{abstract}
Let $K$ be a convex body in ${\bf R}^n$ and $B$ be the Euclidean unit ball in
${\bf R}^n$. We show that
$$\mbox{lim}_{t\rightarrow 0} \frac{|K| -|K_t|}{|B| - |B_t|}=
\frac{as(K)}{as(B)},$$
where $as(K)$ respectively $as(B)$ is the affine surface area of $K$
respectively
$B$ and $\{K_t\}_{t\geq 0}$, $\{B_t\}_{t\geq 0}$ are general families of
convex bodies
constructed from $K$, $B$ satifying certain
conditions. As a corollary we get results obtained in [M-W], [Schm],[S-W]
and[W].

\end{abstract}

\vskip 10mm

The affine surface area $as (K)$ was introduced by Blaschke [B] for
convex bodies in ${\bf R}^3$ with
sufficiently smooth boundary and by Leichtweiss [L1] for convex bodies
in ${\bf R}^n$ with
sufficiently smooth boundary as follows
$$as (K) = \int_{\partial K} \kappa(x)^{\frac{1}{n+1}} d\mu(x),$$
where $\kappa(x)$ is the Gaussian curvature in
$x \in \partial K$ and $\mu$ is the surface measure on $\partial K$.
As it occurs naturally in
many important
questions, so for example in the approximation of convex bodies by polytopes
( see the survey article of Gruber [Gr] and the paper by Sch\"utt [S])
or in a priori estimates for PDEs
[Lu-O], one wanted to have extensions of the affine surface area to arbitrary
convex bodies
in ${\bf R}^n$ without any smoothness assumptions of the boundary.
\par
Such extensions were
given in recent years by Leichtweiss [L2], Lutwak [Lu], Meyer and Werner [M-W],
Schmuckenschl\"{a}ger [Schm], Sch\"{u}tt and Werner
[S-W] and Werner [W].
\par
The extensions of affine surface area to an arbitrary convex body $K$ in
${\bf R}^n$
in [L2], [M-W], [Schm], [S-W] and [W] have a common feature:
\newline
first a specific family $\{ K_t\} _{t\geq 0}$ of convex bodies is constructed.
This family is different in each of the extensions [L2], [M-W], [Schm],
[S-W] and [W]
but of course related to the given convex body $K$.
\newline
Typically the families
$\{K_t\}_{t\geq 0}$ are obtained from $K$ through a ``geometric" construction.
In [L2] respectively [S-W] this geometric construction gives as
$\{K_t\}_{t\geq 0}$
the family of the floating
bodies respectively the convex floating bodies. In [M-W] the geometric
construction
gives the family of the Santal${\mbox{\'o}}$-regions, in [Schm] the
convolution bodies and in
[W] the family of the illumination bodies.
\par
The affine surface area is then obtained by using expressions
involving
volume differences $|K|-|K_t|$ respectively $|K_t|-|K|$.
\vskip 3mm
Therefore it seemed natural to ask whether there are completely general
conditions on a family $\{ K_t\} _{t\geq 0}$ of convex bodies in ${\bf R}^n$
that (in connection with
volume difference expressions) will give us affine surface area.
We give a positive answer to this question which was asked - among others -
by A. Pe\l czy\'{n}ski.
\vskip 10mm

Throughout the paper we shall use the following notations.
\newline
$B(a,r)=B^n(a,r)$ is the n-dimensional Euclidean ball with radius $r$
centered at $a$. We put $B = B(0,1)$.
By $||.||$ we denote the standard
Euclidean  norm on ${\bf R}^n$, by $<^., ^.>$ the standard inner product on
${\bf R}^n$.
For two points $x$ and $y$ in
${\bf R}^n$ $[x,y]=\{\alpha x +(1-\alpha )y: 0 \leq \alpha \leq 1\}$
denotes the
line segment from x to y. For a convex set
$C$ in ${\bf R}^n$ and a point $x \in {\bf R}^n \setminus C$, co$[x,C]$ is the
convex hull of $x$ and $C$.
\newline
$ \mathcal {K}$ denotes the set of convex bodies in ${\bf R}^n$.
For  $K\in {\mathcal K}$,  $\mbox{int}(K)$ is the interior of $K$
and $\partial K$ is the boundary of $K$. For $x \in
\partial K$, $N(x)$ is the outer unit normal vector to $\partial K$ in $x$.
We denote the n-dimensional volume of K by $\mbox{vol}_n(K) = |K|$.
\newline
Let $K\in \mathcal{K}$ and $x \in \partial K$ with unique outer unit normal
vector $N(x)$.
We say that $\partial K$ is approximated in $x$ by a ball from the inside
(respectively from the outside)
if there exists a hyperplane $H$ orthogonal to $N(x)$ such that $H \cap
\mbox{int}(K)
\neq \emptyset$ and a Euclidean ball
$B(r)=B(x-rN(x),r)$ (respectively $B(R)=B(x-RN(x),R)$) such that
$$B(r) \cap H^+ \subseteq K \cap H^+ $$
respectively
$$K \cap H^+ \subseteq B(R) \cap H^+ .$$
Here $H^+$ is one of the two halfspaces determined by $H$.

\newpage
\medskip\noindent {\bf Definition 1}
\vskip 4mm
\begin{it}
For $ t \geq 0$, let $ {\mathcal F}_{t} : {\mathcal K} \rightarrow
{\mathcal K}, \hspace {0.1in}
K \longmapsto {\mathcal F}_t(K)=K_t$, be a map with the following properties
\vskip 4mm
(i) $K_0 =K$ and
\newline
either $K_t \subseteq K$ for all $t \geq 0$ and $ {\mathcal F}_{t}$ is
decreasing in $t$
(that is $K_{t_1} \subseteq  K_{t_2}$ if $t_1 \geq t_2$)
\newline or
$K \subseteq K_t$ for all $t \geq 0$ and  $ {\mathcal F}_{t}$ is increasing
in $t$.

\vskip 4mm
(ii) For all affine transformations A with $\mbox{det}A \neq 0$, for all $t$
$$(A(K))_{|\mbox{det}A|t} = A(K_t).$$
\vskip 4mm
(iii) For all $t \geq 0$, $B_t$ is a Euclidean ball with center 0 and
radius $f_1(t)$ and
$$lim_{t \rightarrow 0}| \frac{|B| -|B_t|}{t^{\frac{2}{n+1}}}| = c,$$
where c is a constant (depending on n only).

\vskip 4mm
(iv) Let $ x\in \partial K$ be approximated from the inside by a ball $B(r)$.
\newline
If $H^+ \cap \partial(K_t) \cap \partial(B(r))_s \neq \emptyset$ for some
$s$ and $t$,
then $s \leq C t$ where $C$ is a constant (depending only on n).

\vskip 4mm
(v) Let $\epsilon > 0$ be given and $ x \in \partial K$ be such that it is
approximated from the inside by a ball $B(\rho - \epsilon)$ and from the
outside by a ball
$B(\rho + \epsilon)$. There exists a hyperplane $H$ orthogonal to $N(x)$ and
$t_0$ such that whenever
$$H^+ \cap \partial(K_t) \cap \partial(B(\rho - \epsilon))_s \neq \emptyset,
\hspace{0.1in} \mbox{for} \hspace{0.1in}t \leq t_0, \hspace{0.1in}s=s(t),$$
respectively
$$H^+ \cap \partial(K_t) \cap \partial(B(\rho + \epsilon))_s \neq \emptyset,
\hspace{0.1in} \mbox{for}\hspace{0.1in} t \leq t_0, \hspace{0.1in}s=s(t),$$
then
$$s \leq (1 + \epsilon) t$$
respectively
$$s \geq (1 - \epsilon) t.$$

\end{it}

\newpage
\medskip\noindent {\bf Remarks 2}
\vskip 4mm
(i) Note that the maps ${\mathcal F}_t$ are essentially determined by the
invariance
property 1 (ii) and by their behaviour with respect to Euclidean balls.

\vskip 4mm
(ii) Let $f_r(t)$ be the radius of $B(0,r)_t$. Then it follows immediately
from
Definition 1 (ii), (iii) that
$$lim_{t \rightarrow 0} \frac{r - f_r(t)}{1 - f_1(t)} = r^{\frac{n-1}{n+1}}.$$

\vskip 4mm
(iii) For some examples the following Definition 1$^\prime$ is  easier to check
than Definiton 1.
\vskip 3mm
{\bf Definition 1$^\prime$}
\newline
\newline
(i) - (iii) as in Definition 1.
\newline
(iv)$^\prime$ If $s < t$, then $K_t \subseteq \mbox{int} (K_s)$.
\newline
(v)$^\prime$ If  $K \subset L$ where $L$ is a convex body in ${\bf R}^n$,
then $K_t \subseteq L_t$
for all $t \geq 0$.
\vskip 4mm
However not all the examples mentioned below satisfy (iv)$^\prime$ and
(v)$^\prime$.
For instance the  illumination bodies (defined below) do not satisfy
(v)$^\prime$.
\vskip 7mm

\medskip\noindent {\bf Examples for Definitions 1 and 1$^\prime$}
\vskip 4mm
1. The (convex) floating bodies [S-W]
\vskip 3mm
Let $K$ be a convex body in ${\bf R}^n$ and $t \geq 0$.
$F_t$ is a (convex) floating body if it is the intersection of all
half-spaces whose defining
hyperplanes cut off a set of volume $t$ of $K$. More precisely, for $u \in
S^{n-1}$
let $a_t^u$ be defined by
$$t=|\{x \in K: <u,x> \geq a_t^u \}|.$$
Then
$$F_t =  \cap _{u \in S^{n-1}} \{x \in K: <u,x> \leq a_t^u \}$$
is a (convex) floating body.
\newline
The family $\{F_t\}_{t \geq 0}$ satiesfies Definitions 1 and 1$^\prime$.

\vskip 5mm
2. The Convolution bodies [K], [Schm]
\vskip 3mm
Let $K$ be a symmetric convex body in ${ \bf R}^n$ and $t \geq 0$. Let
$$C(t) = \{ x \in { \bf R}^n: |K \cap (K+x)| \geq 2t \}$$
and $$C_{t} = \frac{1}{2} C(t).$$
Then $\{C_t\}_{t \geq 0}$ satisfies Definitons 1  and 1$^\prime$.

\vskip 5mm
3. The Santal${\mbox{\'o}}$-regions [M-W]
\vskip 3mm
For $t \in {\bf R}$ and a convex body $K$ in ${\bf R}^n$ the
Santal${\mbox{\'o}}$-region $S(K,t)$
of $K$ is defined as
$$S(K,t) =\{x \in  K: \frac{|K||K^x|}{|B|^2} \leq t\},$$
where $K^x=(K-x)^0 =\{z \in {\bf R}^n: <z, y-x>\hspace {.1in} \leq 1
\hspace {.1in}\mbox{for all } y \in K\}$
is the polar of $K$ with respect to $x$.
(We consider only these t for which $S(K,t) \neq \emptyset$).
\newline
Put
$$S_t = S(K, \frac{|K|}{t|B|^2})=\{x \in K: |K^x| \leq \frac{1}{t} \}.$$
Then  the family $\{S_t\}_{t \geq 0} $, satisfies Definitions 1 and
1$^\prime$.

\vskip 5mm
4. The Illumination bodies [W]
\vskip 3mm
Let $K$ be a convex body in ${ \bf R}^n$ and $t \geq 0$.
The illumination body $I_t$ is the convex body defined as
$$I_t = \{ x \in { \bf R}^n: |\mbox{co} [x,K] \backslash K| \leq t \}.$$
Then the family $\{I_t\}_{t \geq 0}$ satiesfies Definition 1.

\vskip 7mm

\medskip\noindent {\bf Theorem 3}

\par
\begin{it}
Let K be a convex body in ${\bf R}^n$.
For all $t \geq 0$ let $K_t$ respectively $B_t$ be
convex bodies obtained from K respectively B by Definition 1 or 1$^\prime$.
Then
$$\mbox{lim}_{t\rightarrow 0} \frac{|K| -|K_t|}{|B| - |B_t|}=
\frac{as(K)}{as(B)}.$$

\end{it}
\vskip 3mm
\medskip\noindent {\bf Remark}
\newline
Note that
$$as(B) = vol_{n-1}(\partial B) = n|B|.$$

\newpage
\medskip\noindent {\bf Corollary 4}
\vskip 4mm
\begin{it}
(i) [S-W]
\newline
Let K be a convex body in ${\bf R}^n$ and
for $t \geq 0$ let $F_t$ be a floating body. Then
$$lim_{t \rightarrow 0}
c_n \frac{|K| -|F_t|}{t^{\frac{2}{n+1}}} = as(K). $$
where $c_n= 2 \hspace{.1in} (\frac{|B^{n-1}|}{n+1})^{\frac{2}{n+1}}.$

\vskip 3mm
(ii) [Schm]
\newline
Let K be a symmetric convex body in ${\bf R}^n$ and
for $t \geq 0$ let $C_{t}$ be a convolution body. Then
$$lim_{t \rightarrow 0}
c_n \frac{|K| -|C_t|}{t^{\frac{2}{n+1}}} = as(K). $$
where $c_n$ is as in (i).

\vskip 3mm
(iii) [M-W]
\newline
Let K be a convex body in ${\bf R}^n$ and
for $t \geq 0$ let $S_{t}$ be a Santal${\mbox{\'o}}$-region. Then
$$lim_{t \rightarrow 0}
e_n \frac{|K| -|S_{t}|}{t^{\frac{2}{n+1}}} = as(K). $$
where $e_n= \frac{2}{|B|^{\frac{2}{n+1}}}.$

\vskip 3mm
(iii) [W]
\newline
Let K be a convex body in ${\bf R}^n$ and
for $t \geq 0$ let $I_t$ be an  illumination body. Then
$$lim_{t \rightarrow 0}
d_n \frac{|I_t| -|K|}{t^{\frac{2}{n+1}}} = as(K). $$
where $d_n= 2 \hspace{.1in}(\frac{|B^{n-1}|}{n(n+1)})^{\frac{2}{n+1}}.$

\end{it}

\vskip 9mm
For the proof of Theorem 3
we need several Lemmas. The basic idea of the proof is as in [S-W].

\newpage
\medskip\noindent {\bf Lemma 5}
\par

\begin{it}
Let K and L be two convex bodies in ${\bf R}^n$ such that $0 \in
\mbox{int} (L)$ and $L \subseteq K$.
Then
\newline
(i)
$$|K| - |L| = \frac{1}{n} \int_{\partial K}
<x, N(x)>(1-(\frac{||x_L||}{||x||})^n) d\mu(x),$$
where $x_L = [0,x] \cap \partial L$ and $\mu$ is the usual surface measure
on $\partial K$.
\newline
(ii)
$$|K| - |L| = \frac{1}{n} \int_{\partial L}
<x, N(x)>((\frac{||x_K||}{||x||})^n -1) d\mu(x),$$
where $x_K$ is the intersection of the half-line from 0 through $x$  with
$\partial K$
and $\mu$ is the usual surface measure on $\partial L$.

\end{it}
\vskip 3mm
The proof of Lemma 5 is standard.

\vskip 7mm
For $x \in \partial K$ denote by $r(x)$ the radius of the  biggest
Euclidean ball contained in $K$ that touches $\partial K$ at $x$.
More precisely
$$r(x)=\mbox{max}\{r:x\in B(y,r) \subseteq K \hspace{.1in}\mbox{for some y}
\in K\}.$$
\vskip 5mm

\medskip\noindent {\bf Remark}

It was shown in [S-W] that

(i)If $B \subseteq K$, then
$$\mu \{ x \in \partial K: r(x) \geq \beta \} \geq (1-\beta)^{n-1}
\mbox{vol}_{n-1} (\partial K)$$

(ii) $$\int_{\partial K} r(x)^{- \alpha} d \mu (x) < \infty \qquad
\mbox{for all} \quad \alpha ,
\quad 0 \leq \alpha < 1$$

\vskip 5mm

\medskip\noindent {\bf Lemma 6}
\par
\begin{it}
Suppose $0$ is in the interior of $K$.  Then we have for all $x$ with
$r(x) > 0$ and  for all $t \geq 0$
$$0 \leq \frac{ < x,N(x) > \left(1- (
\frac{\| x_t \|}{\| x \|} )^n \right)}{n (|B| - |B_t|)}
\leq g(x),$$ where
$\int_{\partial K} g(x) d\mu(x) < \infty$.
\newline
$x_t=[0,x]\cap \partial K$ if $K_t \subseteq K$.
$x_t$ is the intersection of the half-line from 0 through $x$  with
$\partial K_t$
if $K \subseteq K_t.$
\end{it}
\vskip 5mm

\medskip\noindent {\bf Lemma 7}
\begin{it}
Let $x_t$ be as in Lemma 6. Then
$$\lim_{t \rightarrow 0}
\frac{ < x,N(x) > \left(1- (
\frac{\| x_t \|}{\| x \|} )^n \right)}{n (|B| - |B_t|)}
\mbox{ exists a.e.}$$
and is equal to
\vskip 2mm
(i)  $\frac{\rho(x)^{-\frac{n-1}{n+1}}}{n |B|}$ if the indicatrix of Dupin at
$x \in \partial K$ is an $(n-1)$-dimensional sphere with radius $\sqrt{\rho
(x)}$.
\vskip 2mm
(ii)  $0$, if the indicatrix of Dupin at $x$ is an
elliptic cylinder.
\end{it}
\vskip 5mm
\medskip\noindent {\bf Remark}

(i) $r(x) > 0$ a.e.
[S-W] and the indicatrix of Dupin exists a.e. [L2] and is an ellipsoid or an
elliptic cylinder.
\vskip 2mm
(ii) If the indicatrix is an ellipsoid, we can reduce this case to the case
of a sphere by
an affine transformation with determinant 1 (see  [S-W]).
\vskip 5mm
\medskip\noindent {\bf Proof of Theorem 3}

We may assume that $0$ is in the interior of $K$.  By Lemma 5 and with the
notations of Lemma 6 we have

$$\frac{|K| - |K_t|}{|B| - |B_t|} = \frac{1}{n} \int_{\partial K}
\frac{<x, N(x)>(1-(\frac{||x_t||}{||x||})^n)}{|B| - |B_t|} d\mu(x)$$
By Lemma 6 and the Remark preceding it, the functions under the
integral sign are
bounded uniformly in $t$ by an $L^1$-function and by Lemma 7 they are
converging pointwise a.e.
We apply Lebesgue's convergence theorem.

\vskip 7mm

\medskip\noindent {\bf {Proof of Lemma 6}}
\par
Let $x \in \partial K$ such that $r(x) > 0$.
We consider the proof in the case of Definition 1$^\prime$ and of Definition 1
in the case where $K_t \subseteq K$ for all $t \geq 0$.
The case of Definition 1 where $K \subseteq K_t$ for all $t \geq 0$ is
treated in a similar way.
\vskip 3mm
As
$\| x\| \geq \| x_t \|$, we have for all $t$

\begin{equation}
\frac{1}{n} <x,N(x) > \left(1- ( \frac{ \|
x_t
\|}{\| x \|} )^n  \right)
\leq \hspace{.1in}<\frac{x}{\| x \|}, N(x)> \| x \ - x_t\|
\end{equation}
\newline
Put $r(x) = r$, $x-r(x)N(x) = z$ and $<\frac{x}{\| x \|}, N(x)> =
\mbox{cos} \theta$.
\newline
We can assume that there is an $\alpha > 0$ such that
\begin{equation}
B(0,\alpha) \subseteq K \subseteq
B(0,\frac{1}{\alpha}),
\end{equation}
and hence
$$\mbox{cos} \theta ||x -x_t|| \leq \frac{2}{\alpha}.$$
Let $\epsilon > 0$ be given.
By Remark 2 (ii) there exists $t_1$ such that for all $t \leq t_1$
\begin{equation}
r(1-\frac{1-f_1(t)}{r^{\frac{2n}{n+1}}}(1+ \epsilon)) \leq f_r(t) \leq
r(1-\frac{1-f_1(t)}{r^{\frac{2n}{n+1}}}(1- \epsilon)).
\end{equation}
Let $t_0$ be such that $C t_0 < t_1$.
By Definition 1, (i) $f_1(t)$ is decreasing in t, hence we have for all $t
\geq t_0$
$$f_1(t) \leq f_1(t_0)$$
and thus for all $t \geq t_0$ with (1) and (2)
$$\frac{ < x,N(x) > \left(1- (
\frac{\| x_t \|}{\| x \|} )^n \right)}{n (|B| - |B_t|)} \leq
\frac{2}{\alpha |B| (1-(f_1(t_0))^n)}.$$
Therefore the expression in question is bounded by a constant in this case
and hence integrable. It remains to consider the case when $t < t_0$.

\vskip 3mm
a) Suppose  first that
$$\|x-x_t\|\hspace{.1in} < \hspace{.1in}r \quad \mbox{cos} \theta.$$
For $B(z,r)$ we construct the corresponding inner body
$(B(z,r))_s$ such that $x_t$ is a boundary point of  $(B(z,r))_s$.
By Definition 1 (iii) $(B(z,r))_s$ is a Euclidean ball with center z
and radius $f_r(s)$. As $x_t$ is a boundary point of $(B(z,r))_s$,
\begin{equation}
f_{r}(s) =  r(1-\frac{2||x-x_t||\mbox{cos} \theta}{r}+
\frac{||x-x_t||^2}{r^2})^{1/2}
\leq r(1-\frac{||x-x_t||\mbox{cos} \theta}{2 r})
\end{equation}
The last inequality holds by assumption a).
\newline
So far the arguments are the same for Definiton 1 and Definition 1$^\prime$.
From now on they differ slightly.
\newline
By Definition 1 (iv) $s \leq C t$, hence by monotonicity
$f_r(s) \geq f_r(Ct)$ and thus, as $Ct < t_1$, with (3)
$$f_{r}(Ct) \geq
r(1-(1+ \epsilon ) \frac{1-f_1(Ct)}{ r^{\frac{2n}{n+1}}}),$$
which, using Definition 1 (iii) can be shown to be
\begin{equation}
\geq r(1-(1+ \epsilon ) (C^{\frac{2}{n+1}}+ \epsilon )
\frac{1-f_1(t)}{ r^{\frac{2n}{n+1}}}).
\end{equation}
We get from (4) and (5)
\begin{equation}
1-f_1(t) \geq
\frac{||x-x_t||\mbox{cos} \theta \quad r^{\frac{n-1}{n+1}}}
{2 (1+ \epsilon ) (C^{\frac{2}{n+1}}+ \epsilon )}.
\end{equation}
Observe also that
$$|B|-|B_t| = |B| (1-f_1^n(t)) \geq |B| (1-f_1(t)).$$
This inequality together with (1) and (6) shows that
$$ \frac{ < x,N(x) > \left(1- (
\frac{\| x_t \|}{\| x \|} )^n \right)}{n (|B| - |B_t|)}
\leq \frac{2 (1+ \epsilon ) (C^{\frac{2}{n+1}}+ \epsilon )}{
|B|}r^{-\frac{n-1}{n+1}}.$$
And the latter is integrable by the Remark preceding Lemma 6.
\vskip 3mm
In the case of Definition 1$^\prime$ it follows from (iv)$^\prime$ and
(v)$^\prime$ that $s \leq t$.
For if not, then $s > t$, therefore by (iv)$^\prime$ $(B(z,r))_s \subset
\mbox{int} (B(z,r))_t$
and by (v)$^\prime$ $\mbox{int} (B(z,r))_t \subset \mbox{int} (K_t)$, which
contradicts that
$x_t \in \partial K_t \cap \partial (B(z,r))_s$.
Therefore $f_r(s) \geq f_r(t)$ and thus, as $t < t_1$, with (3)
$$f_{r}(t) \geq
r(1-(1+ \epsilon ) \frac{1-f_1(t)}{ r^{\frac{2n}{n+1}}}).$$
We then conclude as above.

\vskip 3mm
b) Now we consider the case when
$$\|x-x_t\|\hspace{.1in} \geq \hspace{.1in}r \quad \mbox{cos} \theta.$$
We choose  $\alpha$ so small that $x_t \not
\in B(0,\alpha)$.
Let $H$ be the hyperplane through 0 orthogonal to $x$.
Then the spherical cone $C = [x,H \cap B(0,\alpha)]$  is contained in $K$
and $x_t \in C$.
Let $d = \mbox{distance} (x_t, C)$. Then
\begin{equation}
d = ||x-x_t||
\frac{\alpha}{(\alpha^2 + ||x||^2)^{\frac{1}{2}}}.
\end{equation}
Let $w \in [0, x_t]$ such that $||x_t-w|| = \frac{d}{2}$. Let $B(w,R)
\subseteq K$
be the biggest Euclidean ball with center w such that $B(w,R) \subseteq K$.
Then $\partial B(w,R) \cap \partial K \neq \emptyset.$ Moreover $R \geq d$,
which implies that $x_t \in B(w,R)$.
Let $(B(w,R))_s$ be the corresponding inner ball such that $x_t \in \partial
(B(w,R))_s$.
\newline
Now we have to distinguish between Definiton 1 and 1$^\prime$.
\newline
By Definition 1, (iv) $ s \leq Ct$.
By monotonicity $f_R(s) \geq f_R(Ct)$
which, as above, is
$$\geq R(1-(1+ \epsilon ) (C^{\frac{2}{n+1}}+ \epsilon )
\frac{1-f_1(t)}{ R^{\frac{2n}{n+1}}}).$$
As $R \geq d$, the latter is
$$ \geq  d(1-(1+ \epsilon ) (C^{\frac{2}{n+1}}+ \epsilon )
\frac{1-f_1(t)}{ d^{\frac{2n}{n+1}}}).$$
On the other hand by construction $f_R(s) = \frac{d}{2}$. Therefore
$$1-f_1(t) \geq \frac{d^{\frac{2n}{n+1}}}{2 (1+ \epsilon)(C^{\frac{2}{n+1}}
+ \epsilon)}.$$
Note also that (2) implies that $\mbox{cos} \theta \geq \alpha ^2$.
Hence with (1), (2), (7) and assumption b) we get that
$$ \frac{ < x,N(x) > \left(1- (
\frac{\| x_t \|}{\| x \|} )^n \right)}{n (|B| - |B_t|)}
\leq \frac{2 (1+\alpha^4)^{\frac{n}{n+1}}(1+ \epsilon)(C^{\frac{2}{n+1}}+
\epsilon)}
{|B| \alpha^{\frac{6n-2}{n+1}}} \hspace{.1in}
r^{-\frac{n-1}{n+1}}
$$
\vskip 3mm
The case of Definition 1$^\prime$ is treated similarly and the above
inequalities
hold true with $C=1$ and
$C^{\frac{2}{n+1}}+ \epsilon=1$.
\vskip 5mm
\medskip\noindent {\bf {Proof of Lemma 7}}
\par
We again consider the case when $K_t \subseteq K $ for all $t \geq 0$ for
Definition 1.
The case $K \subseteq K_t$ for all $t \geq 0$ for Definition 1 and the
case of Definition 1$^\prime$ are done in a similar way (compare the proof
of Lemma 6).

As in the proof of Lemma 6 we can choose $ \alpha > 0$ such that
$$B(0,\alpha) \subseteq K \subseteq
B(0,\frac{1}{\alpha}).$$
Therefore
\begin{equation}
1 \geq \hspace{.1in}< \frac{x}{\| x \|} ,N(x) >\hspace{.1in} \geq \alpha^2.
\end{equation}
We put again $\mbox{cos} \theta = <\frac{x}{||x||}, N(x)>$.
(1) holds, that is
$$\frac{1}{n} < x,N(x) > \left(1- ( \frac{ \|
x_t\|}{\| x \|} )^n  \right)
\leq \hspace{.1in}<\frac{x}{\| x \|}, N(x)> \| x \ - x_t\|$$
Since $x$ and
$x_t$ are colinear,
$$||x||=||x_t||+||x-x_t||$$
and hence
\begin{equation}
\frac{1}{n} < x,N(x) > \left(1- ( \frac{\| x_t\|}{\|
x
\|} )^n \right)=\frac{1}{n} < x,N(x) > \left(
(1-(1-
\frac{\| x -x_t \|}{\| x \|} )^n \right)$$
$$\geq$$
$$< \frac{x}{\| x \|} ,N(x) > \| x -x_t \| \left( 1 -
k_1 \cdot
\frac{\| x -x_t \|}{\| x \|} \right)
\end{equation}
for some constant $k_1$, if we choose $t$ sufficiently large.
\vskip 3mm
(i) \underline{Case where the indicatrix is an ellipsoid}
\newline
We have seen that then we can assume that the indicatrix is a
Euclidean sphere. Let
$\sqrt{\rho(x)}$ be the radius of this sphere.
We put $\rho(x)=\rho$ and we introduce a coordinate system such that $x =0$ and
$N(x) = (0, \ldots 0, -1)$.  $H_0$ is the tangent hyperplane to $\partial K$ in
$x =0$ and $\{ H_{\alpha}: \alpha \geq 0 \}$ is the family of hyperplanes
parallel to $H_0$
that have non-empty
intersection with $K$ and are of distance $\alpha$ from $H_0$.  For $\alpha
> 0$,
$H_{\alpha}^+$ is the
half-space generated by $H_{\alpha}$ that contains $x =0$.
For $a\in {\bf R}$, let $z_a=(0,\ldots 0,a)$ and $B_a=B(z_a,a)$ be the
Euclidean ball
with center $z_a$ and radius a.
As in [W], for
$\varepsilon > 0$ we can choose $\alpha_0$ so small that for all $\alpha
\leq \alpha_0$
$$B_{\rho-\varepsilon} \cap
H_{\alpha}^+
\subseteq
K \cap H_{\alpha}^+
\subseteq
B_{\rho+\varepsilon}
\cap H_{\alpha}^+.$$
We choose t so small that $x_t \in \mbox{int} (B_{\rho-\varepsilon} \cap
H_{\alpha}^+)
(\subseteq  \mbox{int} (B_{\rho+\varepsilon} \cap
H_{\alpha}^+))$.
For $B_{\rho+\varepsilon}$ we construct the corresponding inner body
$(B_{\rho+\varepsilon})_s$
such that $x_t$ is a boundary point of $(B_{\rho+\varepsilon})_s$.
$(B_{\rho+\varepsilon})_s$ is a Euclidean ball with center
$z_{\rho+\varepsilon}$ and radius
$f_{\rho+\varepsilon}(s)$. We have
$$f_{\rho+\varepsilon}(s) = ((\rho+\varepsilon)^2 + \| x -x_t \|^2 - 2
(\rho+\varepsilon)
\| x -x_t \| \mbox{cos} \theta)^{\frac{1}{2}},$$
$$ \geq (\rho+\varepsilon) (1-\frac{\| x -x_t \|\mbox{cos}
\theta}{\rho+\varepsilon}).$$
Definition 1, (v) implies that $s \geq (1- \varepsilon)t$, hence by
monotonicity
$f_{\rho+\varepsilon}(s) \leq f_{\rho+\varepsilon}((1-\varepsilon)t)$,
which, for $t$ small enough
is  (compare with the proof of Lemma 6)
$$ \leq (\rho+\varepsilon) (1-(1-k_2 \varepsilon)\frac{1-f_1(t)}
{(\rho+\varepsilon)^{\frac{2n}{n+1}}}),$$
where $k_2$ is a constant.
Thus
$$1-f_1(t) \leq \frac{||x-x_t||\mbox{cos} \theta
(\rho+\varepsilon)^{\frac{n-1}{n+1}}}
{1-k_2 \varepsilon}.$$
Note that
$$|B|-|B_t| = |B| (1-f_1^n(t)) \leq n |B| (1-f_1(t)).$$
Therefore with (9)
$$\frac{ < x,N(x) > \left(1- (
\frac{\| x_t \|}{\| x \|} )^n \right)}{n (|B| - |B_t|)}
\geq (1-k_2 \varepsilon)(1-k_1\frac{||x-x_t||}{||x||})
\frac{(\rho+\varepsilon)^{-\frac{n-1}{n+1}}}{n |B| }.$$
This is the lower bound for the expression in question.
\newpage
To get an upper bound we proceed similarily.
For $B_{\rho-\varepsilon}$ we construct the corresponding inner body
$(B_{\rho-\varepsilon})_s$
such that $x_t$ is a boundary point of $(B_{\rho-\varepsilon})_s$.
$(B_{\rho-\varepsilon})_s$ is a Euclidean ball with center
$z_{\rho-\varepsilon}$ and radius
$f_{\rho-\varepsilon}(s)$. We have
$$f_{\rho-\varepsilon}(s) = ((\rho-\varepsilon)^2 + \| x -x_t \|^2 - 2
(\rho-\varepsilon)
\| x -x_t \|\mbox{cos} \theta)^{\frac{1}{2}},$$
$$ \leq (\rho-\varepsilon)(1-\frac{\| x -x_t \|\mbox{cos}
\theta}{\rho-\varepsilon}
(1-\frac{\| x -x_t \|}{2(\rho-\varepsilon)\mbox{cos} \theta})
(1+k_3\frac{\| x -x_t \|\mbox{cos} \theta}{\rho-\varepsilon}
(1-\frac{\| x -x_t \|}{2(\rho-\varepsilon)\mbox{cos} \theta}))),$$
for some constant $k_3$, if $t$ is small enough.
Again by Definiton 1 (v) $s \leq (1+ \varepsilon) t$ and therefore
$f_{\rho-\varepsilon}(s)
\geq f_{\rho-\varepsilon}((1 + \varepsilon ) t)$ which with arguments
similar as before is
$$ \geq
(\rho-\varepsilon) (1-(1+k_4\varepsilon)
\frac{1-f_1(t)}{(\rho-\varepsilon)^{\frac{2n}{n+1}}})$$
with a suitable constant $k_4$. Thus
$$1-f_1(t) \geq$$
\begin{equation}
\frac{||x-x_t||\mbox{cos} \theta}{1+ k_4 \varepsilon } (
1-\frac{||x-x_t||}{2(\rho-\varepsilon)\mbox{cos} \theta})(1+\frac{k_3||x-x_t||
\mbox{cos} \theta}
{\rho-\varepsilon}(1-\frac{||x-x_t||}{2(\rho-\varepsilon)\mbox{cos} \theta}))
(\rho-\varepsilon)^{\frac{n-1}{n+1}}.
\end{equation}
Observe now that
\begin{equation}
|B|-|B_t| = |B| (1-f_1^n(t)) \geq n |B| (1-f_1(t))(1-\frac{n-1}{2}(1-f_1(t))).
\end{equation}
We choose $t$ so small that $1-f_1(t) < \frac{2 \varepsilon}{n-1}$.
This together with (1), (10) and (11) implies that
$$\frac{ < x,N(x) > \left(1- (
\frac{\| x_t \|}{\| x \|} )^n \right)}{n (|B| - |B_t|)}$$
$$\leq$$
$$\frac{1+k_4 \varepsilon}{(1-\varepsilon)(1-\frac{||x-x_t||}
{2(\rho-\varepsilon)\mbox{cos} \theta})(1+k_3 \frac{||x-x_t||\mbox{cos} \theta}
{\rho-\varepsilon}(1-\frac{||x-x_t||}{2(\rho-\varepsilon)\mbox{cos} \theta}))}
\hspace{.1in}
\frac{(\rho-\varepsilon)^{-\frac{n-1}{n+1}}}{n |B| }.$$
Note that $\mbox{cos} \theta \geq \alpha^2$ by (8).
\newline
This finishes the proof of Lemma 7 in the case where the
indicatrix is an ellipsoid.
\newpage
(ii) \underline{Case where the indicatix of Dupin is an elliptic
cylinder}
\newline
Recall that then we have to show that
$$\lim_{t \rightarrow 0}
\frac{ < x,N(x) > \left(1- (
\frac{\| x_t \|}{\| x \|} )^n \right)}{n (|B| - |B_t|)}=0.$$
We can again assume
(see [S-W]) that the indicatrix is a spherical cylinder i.e. the product of a
$k$-dimensional plane
and a $n-k-1$ dimensional Euclidean sphere of radius $\rho$.  We can
moreover assume that
$\rho$ is arbitrarily large (see also [S-W]).

By Lemma 9 of [S-W] we then have for sufficiently small $\alpha$ and some
$\varepsilon > 0$

$$B_{\rho - \varepsilon} \cap H_{\alpha}^+
\subseteq K \cap H_{\alpha}^+.$$
Using similar methods, this implies that
$$\lim_{t \rightarrow 0}
\frac{ < x,N(x) > \left(1- (
\frac{\| x_t \|}{\| x \|} )^n \right)}{n (|B| - |B_t|)}=0.$$

\newpage
\medskip\noindent Elisabeth Werner
\newline
Department of Mathematics
\newline
Case Western Reserve University
\newline
Cleveland, Ohio 44106, U.S.A.
\newline
e-mail: emw2@po.cwru.edu
\newline
and
\newline
Universit\'{e} de Lille 1
\newline
Ufr de Mathematique
\newline
59655 Villeneuve d'Ascq, France

\end{document}